\documentclass[12pt]{article}
\usepackage[centertags]{amsmath}
\usepackage{amsfonts}
\usepackage{amssymb}
\usepackage{amsthm}
\usepackage{amsmath,mathrsfs}
\usepackage{amsmath,amscd}
\usepackage{mathtools}
\usepackage[matrix,arrow,curve,cmtip]{xy}
\title{\textbf{Construction and Properties of the Ground State of Natural Phenomena}}
\author{Renaud Gauthier \footnote{2020 Math. Subj. Class: 68P30, 18N60, 55P42, 18G80, 18E35 . Keywords: stable $\infty$-categories, presentable $\infty$-categories, spectra} \\ \\}
\theoremstyle{definition}
\newtheorem{Net}{Proposition}[subsubsection]
\newtheorem{Chipres}[Net]{Proposition}
\newtheorem{Chinerve}[Net]{Proposition}
\newtheorem{ChinerveRmrk}[Net]{Remark}
\newtheorem{CanIsoSpectra}{Proposition}[subsection]
\newtheorem{CanIsoSpectraRmrk}[CanIsoSpectra]{Remark}
\newtheorem{Sppres}[CanIsoSpectra]{Proposition}
\newtheorem{Spst}[CanIsoSpectra]{Proposition}
\newtheorem{SpLoc}[CanIsoSpectra]{Corollary}
\newtheorem{SpLocRmrk}[CanIsoSpectra]{Remark}
\newtheorem{hSpChiTrCat}{Proposition}[subsection]
\newtheorem{tstr}[hSpChiTrCat]{Proposition}

\newtheorem{ExIRV}{Example}[subsection]
\newtheorem{GrStEx}{Proposition}[subsection]
\newtheorem{LocDec}{Lemma}[section]
\newtheorem{roof}[LocDec]{Construction}

\newcommand{\beq}{\begin{equation}}
\newcommand{\eeq}{\end{equation}}

\newcommand{\rarr}{\rightarrow}

\newcommand{\rlarr}{\rightleftarrows}

\newcommand{\Ob}{\text{Ob\,}}

\newcommand{\xrarr}{\xrightarrow}

\newcommand{\cC}{\mathcal{C}}
\newcommand{\cCop}{\cC^{\op}}
\newcommand{\cD}{\mathcal{D}}
\newcommand{\cE}{\mathcal{E}}
\newcommand{\cF}{\mathcal{F}}
\newcommand{\cG}{\mathcal{G}}

\newcommand{\cK}{\mathcal{K}}
\newcommand{\cM}{\mathcal{M}}

\newcommand{\cP}{\mathcal{P}}

\newcommand{\cS}{\mathcal{S}}

\newcommand{\cX}{\mathcal{X}}

\newcommand{\bR}{\mathbb{R}}

\newcommand{\Cat}{\text{Cat}}

\newcommand{\Fun}{\text{Fun}}

\newcommand{\Hom}{\text{Hom}}

\newcommand{\Mor}{\text{Mor}\,}

\newcommand{\op}{\text{op}}

\newcommand{\Set}{\text{Set}}

\newcommand{\uHom}{\underline{\Hom}}

\newcommand{\Alg}{\text{Alg}}

\newcommand{\Catinf}{\Cat_{\infty}}
\newcommand{\CatDinf}{\Cat^{\Delta}_{\infty}}

\newcommand{\Der}{\mathbb{D}\text{er}}

\newcommand{\dStk}{\text{dSt}(k)}

\newcommand{\del}{\partial}

\newcommand{\Fins}{\cF\text{in}_*}

\newcommand{\oT}{\otimes}

\newcommand{\RuHom}{\bR \uHom}

\newcommand{\SetplusD}{\Set^+_{\Delta}}
\newcommand{\SetD}{\Set_{\Delta}}

\newcommand{\Sp}{\text{Sp}}

\newcommand{\ChiR}{\chi_{\mathbb{R}}}
\newcommand{\Chinat}{\chi_{nat}}
\newcommand{\coCFib}{\text{coCFib}}
\newcommand{\Chisuper}{\chi_{\textit{super}}}
\newcommand{\Catpresinf}{\Cat^{pres}_{\infty}}
\newcommand{\Catstinf}{\Cat^{st}_{\infty}}
\newcommand{\Catstinfpres}{\Cat^{st}_{\infty,pres}}
\newcommand{\Sfin}{\cS^{fin}}
\newcommand{\Sfinst}{\Sfin_*}
\newcommand{\Exc}{\text{Exc}}
\newcommand{\Excst}{\Exc_*}
\newcommand{\LFib}{\text{LFib}}
\newcommand{\RFun}{\text{RFun}}
\newcommand{\TC}{T_{\cC}}
\newcommand{\MTC}{\cM^T(\cC)}
\newcommand{\Cpos}{\cC_{\geq 0}}
\newcommand{\Dec}{\text{Dec}}

\newcommand{\Ztwo}{\mathbb{Z}_2}

\begin{document}
\maketitle
\begin{abstract}
	We construct an $\infty$-category $\cG$ as a model for the Ground State of physical phenomena and we provide properties of its manifestations $\chi = \Fun(\cG,\Catinf)$ in $\Catinf$ as well as of its $\infty$-category of spectra $\Sp(\chi)$. 
\end{abstract}

\section{Introduction}
In a series of papers (\cite{RG}, \cite{RG2}, \cite{RG3}, \cite{RG4}) we have used  derived stacks with simplicial algebras for bases to model natural phenomena. The question presently is how do those algebras form in the first place. More generally, if we consider higher constructs such as $\infty$-categories, $\infty$-operads and the like to regulate natural phenomena, a natural question is, how are those produced? The present paper has for aim to answer this question by constructing a ground state from which higher algebraic objects of interest arise.\\ 

This paper is comprised of two parts, the first one being about the construction of this ground state, the second being about the properties of its manifestations in $\Catinf$.\\

For the sake of motivating the existence of a ground state, we first argue that all concepts involved in describing higher algebraic constructs can be dealt with information theoretically. An organized and thorough analysis of such concepts or notions leads to a layering of their intrinsic definitions, producing what we refer to as iterative relative vistas, a vista consisting in seeing an object as a generic object of some ambient category, within which it can be identified relative to other objects, thereby producing what we call a relative vista. This is done in a minimal fashion, that is one only uses that information which is necessary to describe the object of interest. Beyond that if there remain concepts that need to be defined further, the relative vista picture can be applied to such concepts, thereby resulting in an overall inductive process which we refer to as an iterative relative vista, or $IRV$, which provides a layered deconstruction of a given object. From the perspective of a core object $a$ of interest, $IRV(a)$ can be represented as a recombining tree, all of whose branches - representing vistas -  end in $\Ztwo$, the reason being that as we progress along branches, which we refer to as rays, one introduces concepts that are more and more abstract, until the point dichotomies in a relative setting eventually boil down to considering objects of $\Ztwo$. Having $IRV$s ending in $\Ztwo$ signifies that they provide complete information theoretic deconstructions of their core objects. By the same argument, even if initially rays branch out, precisely because the higher we go in the $IRV$-tree the more abstract the concepts become, branches eventually meet again, hence the recombining feature of such a tree. We denote by $\cG_0'$ the category of $IRV$s for all concepts that come into play when dealing with natural phenomena. We will use later that morphisms between objects of $\cG_0'$ are $IRV$s themselves, thereby making $\cG_0'$ into an internal category. Flipping such a picture of an $IRV$ on its head we have an object that is based at $\Ztwo$, the quintessential starting point for defining concepts, and ends at a given concept. This results in a category $\cG_0 \cong \cG_0'$. Both pictures present information in deconstructed form, whether it be with $IRV$s centered at individual concepts under consideration, or in a  more constructive fashion from the ground up with those inverted $IRV$s. If we regard $\Ztwo$ as a building block of objects of $\cG_0$ under a splicing product $\wedge$, since $\cG_0 \cong \cG_0'$, and $IRV$s are information theoretic deconstructions of concepts, one can therefore regard $\cG_0$ as the categorical ground state of all physical phenomena.\\

From there we work with $\cG = N(\cG_0) \in \Catinf$, which we refer to as the (higher categorical) ground state of natural phenomena. Once such a ground state $\cG$ is constructed we focus on its properties, and more particularly on its manifestations $\chi = \Fun(\cG, \Catinf)$, which we show is a presentable $\infty$-category. Since we work from an algebraic perspective, to study $\chi$ further it is natural to consider its category of spectra $\Sp(\chi)$, and we focus on this $\infty$-category next. We will show among other things that this is a presentable, stable $\infty$-category. Finally, still along the same lines, we are then led to considering its associated triangulated category $h\Sp(\chi)$ and we study some of its properties.

\section*{Notations}
We denote by $\Catpresinf$ the $\infty$-category of presentable $\infty$-categories, and by $\Catstinfpres$ the $\infty$-category of stable, presentable $\infty$-categories. We will use \cite{Lu} and \cite{Lu2} extensively, so for the sake of simplicity we will refer to those references as HTT and HA respectively.

\section{Construction of the ground state}

\subsection{Motivation}
In developing a notion of ground state, we focus on natural phenomena, which arise in the setting of Segal topos (\cite{RG}), and originate from simplicial algebras, or more generally $\infty$-categories and $\infty$-operads. If one seeks a ground state of such a formalism, one aims to construct a base setting from which one can generate such higher algebraic objects. Further, one wants this base to have the universal property that any algebraic construct derives from it, hence it cannot be produced from a more elementary setting itself, in other terms it is fundamental. To construct such a ground state, we first observe that differences between various concepts can be made by way of information (read characterization), a language common to all theories, however disparate they are. Thus we place ourselves at the information theoretic level. This has the advantage of putting different formalisms on a same footing.\\

To even consider a setting that predates algebraic constructs indicates that one is considering deconstructions, and this setting we are aiming for being fundamental, this subsumes that its objects cannot be deconstructed further. Additionally, a feature one would like such a base to have is that its objects can no longer be differentiated, since the existence of differences subsumes having different characteristics, pointing to possible further deconstructions. Keeping those points in mind, we will construct a categorical ground state $\cG_0'$, whose objects are in fully deconstructed form, on which one has a product $\wedge$ called \textbf{splicing} turning $\cG_0'$ into a symmetric monoidal category, with inverted category $\cG_0$ (for lack of a better word), whose objects are all generated from $\Ztwo$, from the perspective of which objects are undifferentiated. All those concepts will be introduced below. For the moment suffices it to say objects of $\cG_0'$ are fully deconstructed concepts, and $\cG_0$ has $ \mathbb{Z}_2$ for basis, from which all concepts can be produced by splicing.

\subsection{Definitions}
\subsubsection{$IRV$s - Construction}
For the sake of deconstructing a given concept, we consider its attributes separately and we map out their inter-relationships in an organized fashion, leading to a structured description, according to the following scheme.\\

In a first time we adopt the position according to which every concept $a$ can be seen as an object of some ambient category $\cC$ which will be referred to as a \textbf{vista}. To have a vista allows one to put $a$ in perspective, to open our vista on that object, hence the name. Observe that the act of seeing $a$ as being embedded in some ambient category $\cC$ suffices to define a vista, without any need to characterize $a$ in $\cC$ other than being a generic object of $\cC$. Note that this is similar in spirit to using a functor of points approach, which is functorial, whereas here one considers a categorical construct. For the sake of having a more precise picture at this level, that is for the purpose of identifying the exact type of $a$ within $\cC$, for characterizing the object of interest within this ambient category, one has to compare $a$ to other objects of $\cC$. Thus we seek a relative point of view. We call a \textbf{relative vista} for a concept $a$ a vista $\cC$ within which $a$ is compared to other objects of $\cC$. It is important to emphasize that in a relative vista we only use objects $b$ that have a connection with $a$, either through a perception $b \rarr a$, or a manifestation $a \rarr b$. Actually one has to be a bit more general as we will see later and we have to allow that we consider any object $b$ with a connection to $a$. For instance if $a$ is a morphism and $b$ is another morphism that composes with $a$, $b$ should figure in the relative vista for $a$, in spite of the fact that there may not be a manifest morphism $b \rarr a$. What is important is that we only keep those objects of a vista with a connection with the core object of interest. Any additional information would not provide us with an optimal information theoretic characterization of $a$. Going back to our construction of vistas for our object $a$, suppose  $b \in \cC$ is another object, and there is a morphism $f:b \rarr a$ of $\cC$ for illustrative purposes that sheds some light on the nature of $a$. Such a morphism possibly introduces additional concepts which become de facto part of the characterization of $a$. Those concepts themselves can be put into light by using their own relative vista. Due to the iterative nature of this process, we arrive at the concept of \textbf{iterative relative vista}, or $IRV$, providing a layered description of a core object. For our purposes we will only consider complete descriptions, that is iterative deconstructions till there is no non-trivial concept left to be defined, and we will reserve the acronym $IRV$ for such complete descriptions. To consider such $IRV$s is essential since our ultimate aim is to fully deconstruct concepts, and stopping midway at some point where we would leave implied concepts undefined would leave us with something that is not fully analyzed. Thus because we aim for an ultimate ground state of all natural phenonena, each given concept involved in defining natural phenomena must be deconstructed till there are no dichotomies involved in its description, in a sense we will make precise below.\\
\begin{ExIRV}
	We provide an elementary example of an $IRV$. Consider the torus $T^2$. It is a manifold, so we see it as an object of the category of manifolds. This provides a vista for $T^2$. But a manifold is first and foremost a topological space. This provides a vista for manifolds. A space in turn is a set with some additional structure. From there one can branch out by considering sets, or structures. Structures are decorations on preexisting objects. Decorations are inherently defined by a before and an after, which can be encoded by $\Ztwo$. So far we have iterative vistas, but nothing is relative yet. Coming back to manifolds now, within the category of manifolds, there are other objects, such as the sphere $S^2$. Relative to $S^2$, $T^2$ is different in that it has genus 1. This partly identifies $T^2$ relative to another object of the category of manifolds. A genus is an integer, object of $\mathbb{Z}_+$, itself a group, itself a set with additional structure, and here we see we have a ray that winds back to decorations, thus branches of an $IRV$ recombine. Focusing on the genus again, it is used in the Euler characteristic, which is a topological invariant, etc...\\ 
\end{ExIRV}

\subsubsection{General considerations}
A word of caution is needed at this point; if we want to fully dissect natural phenomena, we are limited by our understanding of physical laws, so naturally what we would consider to be $IRV$s are only truncations thereof. Therefore in this work we place ourselves from the perspective of a Source which has access to a universe $\cK$ of all knowledge within which $IRV$s can be constructed.\\

Note that one would be tempted to characterize $IRV$s as being graded objects. However if it is true that there is a gradation in the construction of such objects, it is by no means linear since we are really looking at tree-like objects, hence we limit ourselves to just describing such objects as being inductive, hence the name iterative relative vista, as opposed to graded relative vista for example. \\

Observe that starting from a core concept $a$, one has \textbf{rays} emanating from $a$ each of which consists in a path obtained by following vistas and relative vistas and their successive iterations. Naturally this has for aim to deconstruct information. Note that this is gradual; lower strata of an $IRV$ are more structured, higher levels progressively appear more deconstructed as we move away from the core $a$ of $IRV(a)$. We claim that rays ultimately end by $\Ztwo$ (and may even recombine earlier), the reason being that concepts arise from differences, and differences ultimately can be characterized by using elements of $\Ztwo$. For concepts that are characterized by smooth parameters $\lambda \in \bR$, using Dedekind cuts we can still define such concepts using $\Ztwo$. Note that to have trees ending in $\Ztwo$ also allows us to characterize what is meant by complete decompositions. We will only consider complete $IRV$s whose rays ultimately end by $\Ztwo$, since failing that our core objects are not fully analyzed. It is also in that regard that deconstructed concepts $a$ in the form $IRV(a)$ with rays all ending in $\Ztwo$ appear as undifferentiated concepts since from the perspective of $\Ztwo$ one can argue core concepts $a$ have dissolved into elementary dichotomies, at which point there are no distinctions left other than basic, binary dualities of the most basic kind. \\

\subsubsection{Representations of $IRV$s}
Presently we would like to represent $IRV$s. We start with an object $a$. In a first time we can see it as an object of some ambient category $\cC_0$, which we can represent as $a - \cC_0$. Suppose $\cC_0$ itself is an object of a larger ambient category $\cC_1$. Then by iteration we consider $a - \cC_0 - \cC_1$. Each branch shows a transition from one level to the next. What we have here is a two-steps iterative vista with a single branch. To come back to $a - \cC_0$, we should write this in the form $a - (\cC_0)$. We adopt a relative point of view at this point. Say within $\cC_0$ we have a morphism $b \xrarr{f} a$. Then our vista $a - (\cC_0)$ is promoted to the status of relative vista which now reads $a - (\cC_0,b,f)$. More generally if perceptions are denoted $a^- \rarr a$ and manifestations are denoted $a \rarr a^+$, then a relative vista for $a$ can be represented in the form:
\beq
a - (\cC_0,a^-_1, f^-_1, \cdots, a^{-}_N, f^-_N, a^+_1,f^+_1, \cdots, a^+_P, f^+_P) \nonumber
\eeq
Now we can iterate this picture. Say we have a relative vista for some $a_i^{\pm}$ above:
\beq
a^{\pm}_i - (\cD_1,a^{\pm,-}_{i1},f^{\pm,-}_{i1}, \cdots, a^{\pm,-}_{iQ}, f^{\pm,-}_{iQ}, a^{\pm,+}_{i1}, f^{\pm,+}_{i1}, \cdots) \nonumber 
\eeq
Inserting this into the above relative vista for $a$ produces a 3d-tree. As the categories involved become more and more abstract, rays are bound to recombine, giving the picture of a recombining tree with top at $\Ztwo$. Note also that core objects can be collections of objects as in $IRV(a_1 \cup \cdots \cup a_N)$, in which case we are looking at a multi-trunked tree. Finally, in a still more abstract treatment, one has to allow connections of any type $b - a$ in lieu of morphisms $b \rarr a$, and this is in particular pertinent when discussing morphisms. Say $a \xrarr{f} b \xrarr{g} c$ within a category $\cC$. Then one replaces the simple vista $f - (\cC)$ with the relative vista $f - (\cC,g)$ for instance.\\

\subsubsection{Completions}
Regarding $IRV$s, we work within the universe $\cK$ introduced above. Once an $IRV$ is constructed for an object $a$, this is done once and for all so there is no need to repeat this work. We place ourselves in the position where all $IRV$s have already been constructed. Thus we can replace objects with their corresponding $IRV$s. Such an operation we refer to as \textbf{completion} $\Bbbk$, and sometimes for convenience we will write $\Bbbk(a) = \overline{a}$ for $IRV(a)$, and the completion map is denoted $\Bbbk: a \leadsto \overline{a} = IRV(a)$.\\

\subsubsection{Morphisms of $IRV$s}
To summarize, $IRV$s provide minimal, complete deconstructions of concepts. We denote by $\cG_0'$ the category of $IRV$s, with morphisms defined like so: for $a$ and $b$ two concepts, we have a morphism $IRV(a) \rarr IRV(b)$ if there is a morphism from any concept $a_i$ from any ray of $IRV(a)$ to any concept $b_j$ from any ray of $IRV(b)$, and such a morphism is an $IRV$ itself (since we work within $\cK$, where objects are replaced by their completions), thereby making $\cG_0'$ into an internal category. This is of course assuming one can find an ambient category $\cC_{ij}$ within which we have a morphism $f_{ij}: a_i \rarr b_j$. Regarding compositions, consider three objects $a, b$ and $c$ with respective $IRV$s $IRV(a)$, $IRV(b)$ and $IRV(c)$, and say we have morphisms $a_i \xrarr{\psi_{ij}} b_j \xrarr{\phi_{jk}} c_k$ with obvious notations. We will replace those morphisms by their $IRV$s, and for the purpose of composing morphisms we therefore need to discuss composites of $IRV$s. Those (as we will show) are conveniently given using the \textbf{splicing product} $\wedge$, which we now define. In a first time for $\Sigma_1$ and $\Sigma_2$ two information theoretic objects such that:
\begin{align}
	\Sigma_1 & = (\Sigma_0,\Sigma_1') \nonumber \\
	\Sigma_2 & = (\Sigma_0,\Sigma_2') \nonumber
\end{align}
where $\Sigma_0$ are the common objects between $\Sigma_1$ and $\Sigma_2$, which may occur early or later on in their respective decompositions, then:
\beq
\Sigma_1 \wedge \Sigma_2 = (\Sigma_0, \Sigma_1' \coprod \Sigma_2') \nonumber
\eeq
In other terms the splicing of two such objects amounts to gluing along common objects. Observe that the underlying set of $\Sigma_1 \wedge \Sigma_2$ is:
\beq
\underline{\Sigma_1 \wedge \Sigma_2} \cong \Sigma_1 \cup \Sigma_2 \nonumber
\eeq
Now for $IRV$s $\overline{a}$ and $\overline{b}$ that have points in common, rays from those common points or objects cannot branch out in $\overline{a}$ or $\overline{b}$, so we have decompositions of the form:
\begin{align}
	\overline{a} &= (a_0, \Sigma_0) \nonumber \\
	\overline{b} &= (b_0, \Sigma_0) \nonumber
\end{align}
so that:
\beq
	\overline{a} \wedge \overline{b} = (a_0 \coprod b_0, \Sigma_0) \nonumber
\eeq
From this it follows that set-theoretically:
\beq
	\overline{a} \wedge \overline{b} = \overline{a} \cup \overline{b} = \overline{a \cup b} \nonumber
\eeq
a statement that is preserved when we lift this to the level of structured decompositions since splicing does not modify the order along rays.\\

Now going back to morphisms above, it follows from these considerations that:
\beq
	\overline{\psi_{ij}} \wedge \overline{\phi_{jk}} = \overline{\psi_{ij}} \cup \overline{\phi_{jk}} = \overline{ \psi_{ij} \cup \phi_{jk}} \ni \phi_{jk} \circ \psi_{ij} \nonumber
\eeq
Further taking a completion it follows $\overline{\psi_{ij} \cup \phi_{jk}} \supset \overline{\phi_{jk} \circ \psi_{ij}}$. However in deconstructing $\phi_{jk} \circ \psi_{ij}$, one needs to discuss each of $\phi_{jk}$ and $\psi_{ij}$ individually, so:
	\beq
	\overline{\phi_{jk} \circ \psi_{ij}} \supset \overline{\phi_{jk}} \cup \overline{\psi_{ij}} = \overline{\phi_{jk} \cup \psi_{ij}} \nonumber
	\eeq
	Hence:
	\beq
	\overline{\phi_{jk}} \wedge \overline{\psi_{ij}} = \overline{\phi_{jk}\circ \psi_{ij}} \nonumber
	\eeq
We are looking at:
\beq
\xymatrix{
	\overline{\phi_{jk}} \oT \overline{\psi_{ij}} \ar[r]^-{\wedge} & \overline{\phi_{jk}} \wedge \overline{\psi_{ij}} = \overline{\phi_{jk} \circ \psi_{ij}} \\  
\phi_{jk} \oT \psi_{ij} \ar[u]^{\Bbbk \oT \Bbbk} \ar[r]_{\circ} & \phi_{jk} \circ \psi_{ij} \ar[u]_{\Bbbk} 
} \nonumber
\eeq
Thus composition of morphisms in categories translates at the $IRV$-level to splicing. This also shows that composition of morphisms in $\cG_0'$ is well-defined.

\subsubsection{Algebra of $IRV$s}
We claim the splicing $\wedge$ on $\cG_0$ (or $\cG_0'$) determines a symmetric monoidal structure. This product is clearly symmetric with unit $\Ztwo$ by the very definitions of $\wedge$ and $\Ztwo$, with associator $\alpha: (\Sigma \wedge \Sigma_2) \wedge \Sigma_3 \xrarr{\cong} \Sigma_1 \wedge (\Sigma_2 \wedge \Sigma_3)$ evident on the underlying sets:
\beq
\underline{(\Sigma_1 \wedge \Sigma_2) \wedge \Sigma_3} = (\Sigma_1 \cup \Sigma_2) \cup \Sigma_3 = \Sigma_1 \cup (\Sigma_2 \cup \Sigma_3) = \underline{\Sigma_1 \cup (\Sigma_2 \cup \Sigma_3)} \nonumber
\eeq
and splicing preserves the ordering on objects so we can lift this statement to the level of structured objects, so that $\alpha$ is a bona fide isomorphism.
Due to the nature of the objects in $\cG_0'$ and the definition of morphisms therein, one can see that $\cG_0'$ is a higher category. To use the formalism of \cite{RG5}, if $\cC \in \Catinf$, and $\cC_0$ is its deconstruction, one has the following picture:
\beq
\xymatrix{
	\cC \ar[d]_{dec} \\
	\cC_0 \doteq IRV(\cC) \in \cG_0'
} \nonumber
\eeq
where the $dec$ map has for sole purpose to present $\cC$ information theoretically, something that can be achieved by layering its structure, thereby producing a dynamic and constructive decomposition, which we argue can be presented as an $IRV$. We denote by $\cG_0$ the ``inverted" category for $\cG_0'$, whose objects are $IRV$s presented from a bottom up approach, starting at $\Ztwo$ and ending in a structured object $a$, with $\wedge$ making $(\cG_0,\wedge)$ into a symmetric monoidal category as well. $\cG_0'$ presents concepts as $IRV$s, fully analyzed, layered deconstructions. In contrast, in $\cG_0$ objects are built from $\Ztwo$, and their splicings produce various concepts. Those two categories offer dual pictures. In one, provided by $\cG_0'$, going away from $a$ in $IRV(a)$ we abstract the description of the concept $a$ till we arrive at $\Ztwo$, something we can represent by $\del_{\infty} IRV(a) = \Ztwo$, while in $\cG_0$, by inductively splicing elementary objects thereof we specialize towards more refined concepts. $\cG_0$ is probably a better definition for a ground state of natural phenomena than is $\cG_0'$, since its objects are built from individual bits by splicing, and consequently can be regarded as a category generated by undifferentiated concepts, something one would expect is a characteristic feature of a ground state.\\

\subsection{Existence}

\begin{GrStEx}
	There exists a unique categorical ground state for natural phenomena up to isomorphism.
\end{GrStEx}
\begin{proof}
In \cite{RG5} we defined deconstruction maps that are meant to be understood at the information theoretic level. All the various theories of interest for us will occur in $\Catinf$. Suppose we have a morphism $f: \cC \rarr \cC'$ in $\Catinf$ for instance. One can clearly deconstruct all three objects $\cC$, $\cC'$ and $f$ as follows:
\beq
\xymatrix{
	\cC \ar[dr]^{dec} \ar[rr]^f && \cC' \ar[dl]_{dec} \\
	& \cC_0 \cup f_0 \cup \cC'_0  \in \cG_0'
} \nonumber
\eeq
where the zero subscript denotes deconstructed information in the form of $IRV$s, $\cG_0'$ denotes the yet to be defined ground state, and the union is meant to indicate that we collect together bundles of information. One can present this in the following form:
\beq
\xymatrix{
	F(C) = \cC \ar[r]^{f = Ff_0} & \cC'= F(C') \\
	\cG_0' \ni C = \cC_0 \ar[u]^{F = rec} \ar@{.>}[r]_-{f_0} & C' = \cC'_0  \ar[u]_F
} \nonumber
\eeq
	The various packages of deconstructed information such as $\cC_0$, $f_0$ and $\cC'_0$, are each reconstituted into $\cC$, $f$ and $\cC'$ respectively, and the manner in which this is performed is dubbed $F$. If in $\Catinf$ we have other such morphisms involving $\cC$ or $\cC'$, we can link those up down at the level of $\cG_0'$ to extend the definition of $F$. Another reconstitution would yield another ``functor" $F'$. In this manner we define a category $\cG_0'$ that is such only insofar as the initial objects in $\Catinf$ it is generated from form diagrams. Since $\cG_0'$ arises from the deconstruction of diagrams in $\Catinf$, it appears as a higher category itself, or via the nerve functor we can see it as an object of $\Catinf$. The above picture with reconstituting functors $F$ also shows that $\cG_0'$ qualifies as a ground state from which higher algebraic objects are formed.\\

Suppose now there are two such fundamental ground states $\cG_{0,1}'$ and $\cG_{0,2}'$ which we regard as objects of $\Catinf$ with a functor between them. Using the above argument about deconstructions we can find an underlying ground state in $\Catinf$, providing a deconstruction of such a morphism. Suppose finally we have two ground states, with no obvious morphism between them. To say we have two ground states means one can differentiate one from the other, at the information theoretic level. The information content being non empty one can deconstruct such information. Thus one can still formally form a more elementary, underlying ground state. In this manner one ultimately arrives at a categorical ground state, which we denote by $\cG_0'$, and which is unique up to isomorphism, since all complete deconstructions should produce the same objects in the ground state, information theoretically that is. Further we can also regard this category as an $\infty$-category as is usual in higher category theory.
\end{proof}
Observe that in the above picture one can focus on the formation of $IRV$s under deconstruction, so $f_0$ is therefore an $IRV$, which explains why we defined morphisms in $\cG_0'$ to be in deconstructed form, that is, $IRV$s themselves.\\

\subsection{$\infty$-category of deconstructed info}
In what follows we will work with $\cG_0'$ or (equivalently) $\cG_0$, depending on our point of view. Further, $\cG_0$ corresponds better to a notion of ground state within which one finds fundamental bits of information (namely $\Ztwo$) from which higher concepts can be produced. For the sake of studying this ground state within the context of $\infty$-categories, we take the nerve of $\cG_0$, effectively producing an $\infty$-category, which we refer to as the (working) ground state:
\beq
\cG = N(\cG_0) \nonumber
\eeq

\section{Properties}
\subsection{Ground state}
Presently we do not need to have an explicit definition of $\cG_0$ other than it is a category from which one can generate the $\infty$-categories involved in the production of natural phenomena. Since $\cG_0$ is a category, as is usual in higher category theory, we identify $\cG_0$ with its nerve $N(\cG_0)$ as we just did, so saying we investigate the properties of the ground state, this should be understood to be in reference to $\cG_0$, or equivalently $\cG$.\\

In a first time, by the $\infty$-categorical Grothendieck construction of HTT, conveniently summarized in \cite{AMG} at the level of $\infty$-categories, we have:
\beq
\Fun(\cG, \Catinf) \simeq \coCFib(\cG) \nonumber
\eeq
which recall from \cite{RG5} means that the structure of $\cG$ is holographically accessible in $\coCFib(\cG)$. This is an important observation since it indicates that the ground state is peripherally accessible in a pedestrian sense.\\

\subsection{Manifestations}

\subsubsection{General definitions}
\begin{itemize}
	\item \textbf{$IRV$s produce $\infty$-categories.}
By manifestations of $\cG = N(\cG_0)$, we mean the functor category $\Fun(\cG,\Catinf)$. Functors therein are really reconstruction maps strictly speaking, they compress $IRV$s into various objects they are deconstructions of. We argue such rec maps are enough to produce $\infty$-categories. In the proof of the existence of a ground state, we started by considering the deconstruction of a map in $\Catinf$ $\cC \xrarr{f} \cC'$, producing three $IRV$s $\cC_0$, $f_0$ and $\cC'_0$. Presently we adopt the same reasoning and we start with a morphism $f:a \rarr b$ within an ambient $\infty$-category $\cC$, which we will show can be produced from an $IRV$. Deconstructing $f:a \rarr b$, we get $IRV$s for $a$, $b$ and $f$. But observe that one can splice those $IRV$'s, one based at $a$, the other based at $b$, via the $IRV$ for $f$, a morphism of $\cC$. This means, identifying unions with splicings by abuse of notation:
		\begin{align}
			IRV(a) \cup IRV(f) \cup IRV(b) &\cong IRV(a) \wedge_{\cC} IRV(b) \nonumber \\
			& = \overline{a} \wedge \overline{b} \nonumber \\
			&= \overline{a \cup b} \in \cG_0 \nonumber
		\end{align}
iterating this picture, the deconstruction of an $\infty$-category is isomorphic to an $IRV$, so $IRV$s can functorially map to $\infty$-categories by reconstruction.\\

\item \textbf{Diagrams for natural phenomena.}
We are ultimately interested in the $\cG$-diagrams $G$ in $\Catinf$ that correspond to those algebraic objects that generate natural phenomena. To obtain a well-defined cohesive picture for such $\cG$-diagrams, consider the smallest subcategory of the homotopy category $h\chi$ of $\chi = \Fun(\cG, \Catinf)$ spanned by such functors $G$, and consider its associated $\infty$-category $\ChiR$ following the construction of HTT 1.2.11.\\

Extracting $\ChiR$ from $\chi$ we refer to as \textbf{symmetry breaking}, since we go from $\chi$, at which level all concepts are present, including those that don't intervene in natural phenomena, to an $\infty$-category $\ChiR$ that involves only those diagrams that produce real phenomena.\\

Regarding $\ChiR$ itself, it is important to keep in mind that it contains everything that is necessary to develop a theory of natural phenomena, algebraically speaking that is. Recall that in \cite{RG4} we studied the dynamics of $\cX=\dStk$, the Segal topos that is meant to model natural phenomena, which can be considered to be algebraic. If we consider its $\infty$-categorical instantiation $\Chinat$, it must be generated by $\ChiR$. Indeed, $\cG$-diagrams in $\ChiR$ produce higher algebraic objects involved in natural phenomena by definition, so there are $\cG$-diagrams that specifically produce $\Chinat$. As a matter of fact $\ChiR$ produces quite a bit more, thereby pointing to the fact that $\ChiR$ is actually quite large.\\

Indeed, for the sake of completeness we should also consider an overarching $\infty$-category to $\cX$ above, following \cite{RG4}. Recall that $\cX$ can be identified with $id_{\cX} \in \RuHom(\cX,\cX) = \delta \cX$, and the latter Segal topos corresponds to algebraic deformations of $\cX = \dStk$, thus $\delta \cX$ is also produced by $\ChiR$, or rather its $\infty$-categorical instantiation $\Chisuper$ is, in such a manner that we have the following convenient summarizing diagram:
\beq
\xymatrix{
	& \Chisuper \\
\Chinat \ar@{^{(}->}[ur]  \\
\ChiR  \ar[u]^{rec} \ar@{-}[r] & \ChiR \ar[uu]_{rec} 
} \nonumber
\eeq
		where by $\ChiR \xrarr{rec} \Chinat$ for instance, we have a compact way of saying there are functors of $\ChiR$ that contribute to producing $\Chinat$.\\

		\item \textbf{Dual theories}
To come back to $\ChiR$, among functors $G: \cG \rarr \Catinf$ of interest, some deal with dual theories. Say $A$ and $B$ in $\Catinf$ are dual theories, as well as $A'$ and $B'$, and suppose we have a commutative diagram linking such dual theories as in:
\beq
\xymatrix{
	A \ar@/^/[d] \ar[r] & A' \ar@/^/[d] \\
	B \ar@/^/[u] \ar[r] & B' \ar@/^/[u]
}  \label{CDdual}
\eeq

By definition of $\ChiR$, one can find appropriate functors $F: \cG \rarr \Catinf$ in such a manner that there exists $C \in \cG$, $FC = (A,B)$ and $FC' = (A',B')$ and a morphism $u: C \rarr C'$, so that we have $Fu:FC \rarr FC'$, which further can be decomposed as $Fu_1: A \rarr A'$ and $Fu_2: B \rarr B'$, reproducing the horizontal arrows in \eqref{CDdual}.\\

Focusing on individual dual theories, recall from \cite{RG5} that one can represent dual theories as follows:
\beq
\xymatrix{
	A \ar@/_/[dr]^{dec} \ar@/^/[rr]^{E_{AB}} && B \ar@/_/[ll] \ar@/^/[dl]_{dec} \\
	& C_{AB} \ar@/^1pc/[ul]^{rec}  \ar@/_1pc/[ur]_{rec}
} \nonumber
\eeq
		In our present setting, the reconstruction maps $C_{AB} \rarr A$ and $C_{AB} \rarr B$ above can be repackaged as $C_{AB} = \big (IRV(A) \xrarr{IRV(E_{AB})} IRV(B) \big ) \cong IRV(A) \wedge_{E_{AB}} IRV(B)  = \overline{A \cup B} \mapsto (A \xrarr{E_{AB}} B)$.\\
\end{itemize}

\subsubsection{Properties of $\chi$}

Since $\cG$ is a ground state, and $\chi$ corresponds to a dynamic presentation of its realizations in $\Catinf$, hence is a dual of $\cG$, we would like to argue we have a net neutrality rule, which essentially states that $\sum_{F \in \chi} F$ reproduces $\cG$ in some sense. We make this precise:
\begin{Net}
	There exists a left localization:
	\beq
	\chi \xrarr{L} \cG \nonumber
	\eeq
\end{Net}
i.e. $\cG$ can be recovered from $\chi$ by localization.
\begin{proof}
	From \cite{AMG}, for $\cC \in \Catinf$, we have a commutative diagram in $\Catinf$ with horizontal equivalences:
	\beq
	\xymatrix{
	\Fun(\cC,\Catinf) \ar[r]^{\simeq} & \coCFib(\cC) \\
	\Fun(\cC,\cS) \ar[u] \ar[r]_{\simeq} & \LFib(\cC) \ar[u]
	} \nonumber
	\eeq
	with $\cS$ the $\infty$-category of spaces. Further by Remark 1.2 of the same reference, we have left localization adjunctions: $U: \coCFib(\cC) \rlarr \LFib(\cC): R$. Additionally, by HTT 5.1.3.1, if $K \in \SetD$, $j: K \rarr \cP(K) = \Fun(K^{\op},\cS)$ is fully faithful ($\infty$-categorical Yoneda lemma). Thus applying this to $K = N(\cG_0) = \cG$, we can dynamically represent $\cG$ as $\cP(\cG)$. Collecting things together we have our desired localization:
\beq
	\chi \simeq \coCFib(\cG) \xrarr{U} \LFib(\cG) \simeq \cP(\cG^{\op}) \nonumber
\eeq
\end{proof}

Our first result about $\chi$ proper follows. As an intermediate step we will use that $\Catinf \in \Catpresinf$. One way to show this goes as follows. By HTT 3.1.3.7, if $\SetplusD$ denotes the category of marked simplicial sets, $S \in \SetD$, one can put the cartesian model structure on $(\SetplusD)_{/S}$, which is a combinatorial model structure. By HTT 3.1.4.4 we can turn $(\SetplusD)_{/S}$ into a simplicial category, thereby making $(\SetplusD)_{/S}$ into a simplicial model category. Then use the fact that $(\SetplusD)^{\circ} = \CatDinf$, the simplicial category of small $\infty$-categories, and by definition its nerve $N(\CatDinf) = \Catinf$, that is $\Catinf = N(\SetplusD)^{\circ}$. Finally by HTT A.3.7.6 since $\SetplusD$ is a combinatorial, simplicial model category, it follows $\Catinf$ is presentable.
\begin{Chipres}
	$\chi$ is a presentable $\infty$-category.
\end{Chipres}
\begin{proof}
	It suffices to use the fact that $\Catinf \in \Catpresinf$, $N(\cG_0) \in \SetD$, and HTT 5.5.3.6 which states that if $K \in \SetD$, $\cC \in \Catpresinf$, then $\Fun(K,\cC) \in \Catpresinf$ as well.
\end{proof}
\begin{Chinerve}
	$\chi \simeq N(A^{\circ})$ for some combinatorial, simplicial model category $A$.
\end{Chinerve}
\begin{proof}
	This follows from HTT A.3.7.6 since $\chi \in \Catpresinf$.
\end{proof}
\begin{ChinerveRmrk}
What this says is that $\chi$ is a deconstruction of some homotopical object, or vice-versa, $A$ provides a more compact presentation of $\chi$. $\cG$ being our ground state, $\chi$ corresponding to its functorial presentation in $\Catinf$, objects thereof are reconstructions in the sense of \cite{RG5}. One would be tempted to argue $A$ provides the next step in reconstructions from $\cG$ insofar as it provides a pre-arranged setting. This however is just a ``parallel" reconstruction in a sense, as opposed to functors $G: \cG \rarr \Catinf$ which provide ``horizontal" reconstructions. 
\end{ChinerveRmrk}

\subsection{Spectra}
We consider manifestations of $\cG = N(\cG_0)$ into $\Catinf$, which are packaged into $\chi = \Fun(\cG,\Catinf)$, which is algebraic in nature. If we place ourselves in a higher algebraic context, it becomes natural to consider resolutions and towers of objects at times, and that motivates the appearance of spectra. Thus if we wish to project $\chi$ into higher algebra and study its linearizations, spectra thereof become relevant.\\

By HA 1.4.2.8, for $\cC \in \Catinf$ with finite limits, one can define spectra of $\cC$. Those are reduced and excisive functors $F: \Sfinst \rarr \cC$. Here $\Sfinst$ is the $\infty$-category of pointed objects of $\Sfin$, the smallest full subcategory of $\cS$ stable under finite colimits and containing the final object $*$. Recall also that being excisive means mapping pushout squares to pullback squares, and being reduced means preserving the final object of the base category. We denote by $\Sp(\cC) = \Excst(\Sfinst, \cC)$ the full subcategory of $\Fun(\Sfinst, \cC)$ spanned by spectra of $\cC$.
\begin{CanIsoSpectra}
	We have a canonical isomorphism:
	\beq
	\Sp \chi \cong \Fun(\cG,\Sp(\Catinf)) \nonumber
	\eeq
\end{CanIsoSpectra}
\begin{proof}
This is a direct consequence of HA 1.4.2.9, which states that if $\cC \in \Catinf$ admits finite limits and $K \in \SetD$, then $\Sp(\Fun(K,\cC)) \cong \Fun(K,\Sp(\cC))$. We use this result with $K = N(\cG_0) \in \SetD$, and $\Catinf \in \Catpresinf$, in particular it has finite limits.
\end{proof}
\begin{CanIsoSpectraRmrk}
An informal way to understand what this result says goes as follows. The ground state $\cG$ is presumably that from which algebraic objects involved in natural phenomena originate. $\Sp\chi = \Excst(\Sfinst, \chi)$ being a perception of $\chi$ from the perspective of $\Sfinst$ (and more to the point, is an $\infty$-category as the next result shows), it must be generated from $\cG$. If we just used the fact that $\Sp \chi \in \Catinf$, the only thing we would have is that there exists $G \in \chi$, $C \in \cG$, such that $GC = \Sp \chi$. This however is a punctual representation lacking in inner structure. The above result shows that, owing to the nature of $\Sp \chi$, one can get something more structured, namely that $\Sp \chi$ is actually related to a manifestation $\Fun(\cG,\cE)$ of $\cG$ for some category $\cE$, and that this category is none other than $\Sp(\Catinf)$.\\
\end{CanIsoSpectraRmrk}
\begin{Sppres}
	$\Sp \chi$ is presentable
\end{Sppres}
\begin{proof}
	By HTT 7.2.2.9, if $\cC \in \Catinf$ has a final object, then the $\infty$-category of its pointed objects has a final object, so it is pointed as well. Since $\Sfin$ has a final object, it follows that $\Sfinst = (\Sfin)_*$ is pointed. It also admits finite colimits. Recall also that $\Catinf \in \Catpresinf$. Using these two facts, we can use HA 1.4.2.4 which states that if $\cC$ is a small pointed $\infty$-category admitting finite colimits, $\cD \in \Catpresinf$, then $\Excst(\cC,\cD) \in \Catpresinf$. It follows that $\Sp(\Catinf) = \Excst(\Sfinst, \Catinf) \in \Catpresinf$. By HTT 5.5.3.6, $\Fun(\cG,\Sp(\Catinf))$ is therefore presentable. We conclude by invoking the previous result.
\end{proof}

\begin{Spst}
$\Sp \chi$ is stable.
\end{Spst}
\begin{proof}
We have shown $\chi \in \Catpresinf$, in particular it has finite limits. It follows by HA 1.4.2.17 that $\Sp \chi$ is stable.
\end{proof}

\begin{SpLoc}
There is a small $\infty$-category $\cE$ such that $\Sp \chi$ is equivalent to an accessible left exact localization of $\Fun(\cE,\Sp)$.
\end{SpLoc}
\begin{proof}
Recall that $\Sp = \Sp (\cS_*)$, the $\infty$-category of spectra. The result follows from HA 1.4.4.9 which ensures there is such an $\infty$-category $\cE$, provided $\Sp \chi \in \Catstinfpres$, which is what we have just shown above.
\end{proof}
\begin{SpLocRmrk}
	Another way of stating this is that $\Sp \chi \cong \Fun(\cG, \Sp(\Catinf))$ can be presented by $\Fun(\cE,\Sp)$. This parallelism provides an alternate way of dealing with $\Sp(\chi)$. One can regard $\cE$ as the appropriate translate of $\cG$ so that functors $\cE \rarr \Sp$ reproduce those of the form $\cG \rarr \Sp(\Catinf)$ that form $\chi$. From a relative point of view $\Sp \chi \cong \Fun(\cG, \Sp(\Catinf))$ provides emanations from the ground state, taking $\cG$ as point of reference. If one focuses on $\Sp$ instead, taking it as being a central object, then one can represent $\Sp \chi$ as $\Fun(\cE,\Sp)$. Those therefore provide dual ways of representing the same object, namely $\Sp \chi$, depending on our relative point of view, whether it be that of the ground state, or that of spectra. 
\end{SpLocRmrk}

\subsection{Triangulated categories}
In the continuation of the above algebraic considerations, the next step consists in working with triangulated categories.\\

\begin{hSpChiTrCat}
The homotopy category $h \Sp \chi$ is a triangulated category.
\end{hSpChiTrCat}
\begin{proof}
This follows from HA 1.1.2.15, which states that if $\cC$ is stable, its homotopy category $h\cC$ is a triangulated category. We have just shown that $\Sp \chi$ is stable. The result follows. 
\end{proof}
\begin{tstr}
$\Sp(\chi)$ has a t-structure.
\end{tstr}
\begin{proof}
	Recall from HA 1.2.1.4 that to say $\cC \in \Catstinf$ has a t-structure means its associated triangulated category $h\cC$ has a t-structure. By HA 1.4.4.11 a stable, presentable $\infty$-category can be endowed with a t-structure. Thus $\Sp(\chi) \in \Catstinfpres$ has a t-structure. 
\end{proof}
This will become useful in what follows. As pointed out, spectra provide algebraic decompositions, and HA 7.4.1.23 tells us that if one focuses on certain algebraic aspects of $\Sp(\chi)$, then it appears those are differential in nature. To be precise, this theorem of HA states that $\Fun(\Delta^1,\Alg(\Sp(\chi)))$ is equivalent to an $\infty$-category of algebraic derivations. That this functor category produces a tangent space of some sort comes as no surprise. But that it does correspond to derivations is new information. This shows that spectra provide algebraic resolutions of objects that are differential in character.\\

To use this result we need $\cC = \Sp(\chi)$ to satisfy the conditions of HA 7.4.1.18, that is it needs to be a stable presentable $\infty$-category with a t-structure and a monoidal structure. The latter condition is the only thing we do not currently have. However by HA 1.1.3.4, $\Sp(\chi)$ is stable, so it admits finite products, which means $\Sp(\chi)^{\times} \rarr N(\Fins)$ is a symmetric monoidal $\infty$-category by HA 2.4.1.5, which is further cartesian, that is the associated functor $\pi: \Sp(\chi)^{\times} \rarr \Sp(\chi)$ restricts to an equivalence of $\infty$-categories, so this effectively puts a symmetric monoidal structure on $\Sp(\chi)$. More generally, if $\Catinf^{\oT}$ denotes the $\infty$-category of symmetric monoidal $\infty$-categories, $\Catinf^{\oT,\times}$ the full subcategory spanned by cartesian symmetric monoidal $\infty$-categories, and if $\Catinf^{cart} \subset \Catinf$ denotes the subcategory spanned by $\infty$-categories that admit finite products, then by HA 2.4.1.9 the functor $\theta: \Catinf^{\oT,\times} \rarr \Catinf^{cart}$ is an equivalence of $\infty$-categories, which provides a global statement justifying that one has a well-defined symmetric monoidal structure on $\Sp(\chi)$ as desired.\\

In HA 7.4.1.23 mention is being made of $\Alg(\cC)$. If one wants to focus on the algebraic properties of spectra, one possibility is to consider $\Alg(\cC)$, the category of associative algebra objects of $\cC$. This however is inert. For the sake of introducing some dynamics in the picture one considers $\Fun(\Delta^1,\Alg(\cC)$. HA 7.4.1.23 actually deals with a certain subcategory thereof. Separately the $\infty$-category of derivations is derived from a notion of tangent correspondence, which itself necessitates notions of tangent and cotangent functors. We briefly remind the reader what those are, without getting into too much details.\\

With regards to the tangent bundle, one first needs the notion of stable envelope. By HA 7.3.1.1, a stable envelope of a presentable $\infty$-category $\cC$ is a categorical fibration $U: \cC' \rarr \cC$ with $\cC' \in \Catstinfpres$, $U$ with a left adjoint, and such that for any pointed $\cD \in \Catstinfpres$, one has an equivalence in $\Catinf$: $U_*: \RFun(\cD,\cC') \xrarr{\simeq} \RFun(\cD,\cC)$, which really means that $\cC'$ is a stable adjoint to $\cC$ with $U$ local with respect to pointed stable presentable $\infty$-categories, so effectively one gets a stable resolution of $\cC$. More generally, adopting a relative point of view, if $p: \cC \rarr \cD$ is a presentable fibration, a stable envelope of $p$ is just a pointwise stable envelope in the following sense: one has a categorical fibration $U:\cC' \rarr \cC$ such that $p \circ U$ is still a presentable fibration to preserve the nature of $p$, $U$ carries $(p \circ U)$-cartesian morphisms of $\cC'$ to $p$-cartesian morphisms in $\cC$ for the sake of coherence, and for any $D \in \Ob(\cC)$, $\cC'_D \rarr \cC_D$ is a stable envelope. Armed with this notion one can now define tangent bundles. By HA 7.3.1.9, for $\cC \in \Catpresinf$, a tangent bundle to $\cC$ is a functor $\TC \rarr \Fun(\Delta^1,\cC)$ making $\TC$ the stable envelope of the presentable fibration $p:\Fun(\Delta^1,\cC) \rarr \Fun(\{1\},\cC) \cong \cC$. In other terms, $\TC$ is a fiberwise, stable resolution of the dynamic presentation $\Fun(\Delta^1,\cC)$ of $\cC$. Observe in passing that $\TC = \Exc(\Sfinst,\cC)$.\\

Having the tangent functor, one can now define the cotangent functor. By HA 7.3.2.14, the cotangent complex functor $L: \cC \rarr \TC$ is the left adjoint to $\TC \rarr \Fun(\Delta^1,\cC)$ whose exact construction is given by the definition itself and HA 7.3.2.6. The reader is referred to those for details. \\

The derivations mentioned above combine tangent and cotangent objects by way of a correspondence that we now define. Having $\cC$ and $\TC$, as well as a notion of cotangent complex, if $L:\cC \rarr \TC$ is the associated cotangent complex functor, one can construct, using HA 7.3.6.3, a cocartesian fibration $p:\cM \rarr \Delta^1$ in such a manner that $\cM \times_{\Delta^1}\{0\} \simeq \cC$ and $\cM \times_{\Delta^1} \{1\} \simeq \TC$, with associated functor $\cC \rarr \TC$ being $L$. One denotes $\cM$ by $\MTC$ and refer to it as the tangent correspondence of $\cC$ (see HA 7.3.6.9 for details). Thus $\MTC$ captures $\TC$ and $L$. Then by HA 7.4.1.1, for $\cC \in \Catpresinf$, $p: \MTC \rarr \Delta^1 \times \cC$ a tangent correspondence of $\cC$, a derivation in $\cC$ is a map $f: \Delta^1 \rarr \MTC$ producing a morphism $\eta: A \rarr M$ in $\MTC$ referred to as a derivation of $A$ into $M$. One lets $\Der(\cC) = \Fun(\Delta^1, \MTC) \times_{\Fun(\Delta^1,\Delta^1 \times \cC)} \cC$ be the $\infty$-category of derivations in $\cC$. Exact details defining $f$ above are spelled out in HA. What is important is the nature of those derivations, those are algebraic. The equivalence of HA 7.4.1.23 actually uses a certain subcategory of $\Fun(\Delta^1,\Alg(\cC))$ and a certain type of derivations in $\Der(\Alg(\cC))$, and asks that $\cC = \Sp(\chi)$ satisfies the conditions of HA 7.4.1.18. The reader is referred to those for details. The only thing we need to apply this theorem to our setting is that $\Sp(\chi)$ be monoidal, and this we discussed above. The other conditions of HA 7.4.1.18 are also relatively easy to implement; one asks that the unit object be in $\Cpos$, the tensor product preserves small colimits in each variable, and that it carries $\Cpos \times \Cpos$ into $\Cpos$, and this is something we will assume we have.

\section{Spontaneous reconstruction}
Of independent interest we discuss spontaneous reconstructions. Objects of $\Fun(\cG, \Catinf)$ reconstruct algebraic concepts in $\Catinf$ from $IRV$s in $\cG = N(\cG_0)$, and concepts themselves deconstruct into $IRV$s, two processes which seem inverse to each other. Thus one would like a picture whereby one type of map, be it a reconstruction or a deconstruction, can be inverted. This is what we consider in the present section.\\

Regardless of the situation we are looking at, we will use the same method, namely we will consider a category $\cC$ along with a class $\Dec \subset \Mor \cC$ of deconstruction maps, and we will consider its localization $\cC[\Dec^{-1}]$. For this to be well-defined we want $\Dec$ to be a localizing class. Observe that the choice of dec maps as opposed to rec maps originates from their natural definition.\\

Recall what a localizing class of morphisms consists of (\cite{GM}): in a category $\cC$, with $S$ a class of morphisms in $\cC$, we say $S$ is a localizing class of morphisms if it contains identities, is stable by composition, and has the following extension properties: whenever $s \in S$ and $f \in \Mor \cC$, one can find morphisms $r \in S$ and $g \in \Mor \cC$ that make the following diagrams commutative:
\beq
\xymatrix{
	W \ar@{.>}[d]_r \ar@{.>}[r]^g & Z \ar[d]^s \\
	X \ar[r]_f & Y
} \label{star}
\eeq
and:
\beq
\xymatrix{
	W &Z \ar@{.>}[l]_g \\
	X \ar@{.>}[u]^r & Y\ar[l]^f \ar[u]_s
} \label{2stars}
\eeq
and the last condition is that if $f,g \in \Mor \cC$, if $sf = sg$ for $s \in S$, then $ft = gt$ for some $t \in S$. Observe that we do not ask that this be an equivalence, as is typically the case, because we do not work with the opposite category $\cCop$ of interest, and we focus only on ``left" roofs, and not on right roofs, so this condition is sufficient for our purposes.\\

\begin{LocDec}
If $\cC$ is a category such that $\Mor \cC$ contains a subclass $S$ of deconstruction maps, then $S$ is localizing
\end{LocDec}
\begin{proof}
	The identity is a trivial deconstruction, and since we consider partial deconstructions, their composites are clearly dec maps. The extension condition \eqref{star} is met by observing that deconstructed objects necessarily come from whole objects (or condensates), and so it is for their morphisms, they are deconstructions of morphisms at the level of condensates, hence the lift. For the projection \eqref{2stars}, it just suffices to consider deconstructed maps. For the last condition if two maps $f$ and $g$ are different, $f,g:\cC \rarr \cC$, then there is some $c \in \cC$, $fc \neq gc$. Applying $s \in S$ thereafter as in $sf = sg$ means deconstructing $\cC$, those substructures $fc$ and $gc$ therein that were different are still different in $\cC_0$, the deconstructed form of $\cC$, hence a statement such as $sf = sg$ makes sense only if $f = g$, rendering that condition trivial.
\end{proof}

\begin{roof}
	If $\cC$ is a category with a class of deconstruction maps $S \in \Mor \cC$, then one can define the localization of $\cC$ by $S$ as the category $\cC[S^{-1}]$ of equivalence classes of roofs
	\beq
	\xymatrix{
		&B \ar[dl]_s \ar[dr]^f \\
		A \ar@{-->}[rr] && C
	} \nonumber
	\eeq
where $s \in S$ and $f \in \Mor \cC$, the equivalence between two such ``left" roofs being defined by having a third random overarching roof as in:
\beq
	\xymatrix{
	&& F \ar[dl]_t \ar[dr]^h \\
	&B \ar[dl]_s \ar[drrr]^f && D \ar[dlll]_r \ar[dr]^g \\
	A &&&& E
	}\nonumber
\eeq
making the above figure commutative, with $st \in S$, making $(s,f) \sim (r,g)$. Composition of left roofs $(s,f)$ and $(r,g)$ is provided by the existence of a lift to an overarching roof $(t,h)$ in such a manner that we obtain a third left roof $(st,gh)$ as in:
\beq
	\xymatrix{
	&& F \ar[dl]_t \ar[dr]^h \\
	&B \ar[dl]_s \ar[dr]^f && D \ar[dl]_r \ar[dr]^g \\
	A && C && E
	}\nonumber
\eeq

\end{roof}
We now apply this to $\chi$ objectwise for all functors $F: N(\cG_0) \rarr \Catinf$. Let $\cC = N(\cG_0) \coprod \Catinf$, $\Dec = \Mor(\Ob \Catinf, \Ob N(\cG_0))$, consider $\cC[\Dec^{-1}]$, this generates $\chi$. The choice of dec maps as opposed to rec maps stems from the fact that one starts from a concept, and once it is formed, it can be characterized, hence analyzed, and the act of analyzing a concept is a natural deconstruction in itself. This picture offers one with a spontaneous reconstruction, or equivalently said, fully formed concepts and their corresponding deconstructions coexist simultaneously.

\bigskip
\footnotesize
\noindent
\textit{e-mail address}: \texttt{rg.mathematics@gmail.com}.

\end{document}